\documentclass{amsart}
\usepackage{amsmath}
\usepackage{amssymb}

\usepackage{graphicx}

\newtheorem{rem}{Remark}

\def\am{\mathop{\mathrm{am}}\nolimits}
\def\sn{\mathop{\mathrm{sn}}\nolimits}
\def\cn{\mathop{\mathrm{cn}}\nolimits}
\def\dn{\mathop{\mathrm{dn}}\nolimits}
\def\sd{\mathop{\mathrm{sd}}\nolimits}

\usepackage[pdfauthor={Richard J. Mathar}, colorlinks=true,pdfkeywords={Cone, Sphere, Intersection, Volume},pdftitle={Volume of Intersecting Cone and Sphere},bookmarksopen=true,pdfpagelayout=OneColumn]{hyperref}

\begin{document}

\title{Volume of Intersection of a Cone with a Sphere}

\author{Richard J. Mathar}
\urladdr{https://www.mpia.de/~mathar}
\address{Max-Planck Institute of Astronomy, K\"onigstuhl 17, 69117 Heidelberg, Germany}

\date{\today}

\begin{abstract}
The manuscript provides formulas for the volume of a body
defined by the intersection of a solid cone and a solid sphere
as a function of the sphere radius, of the distance between cone apex and sphere center,
and of the cone aperture angle.

If the sphere center lies on the (extended) cone axis
the analysis may be based on cylinder coordinates fixed
at the cone axis,
and the volume is the sum of the well-known volumes of finite
cones and sphere caps.

At the general geometry the sphere center is \emph{not} on the (extended) cone axis.
Our approach calculates 
the volume by slicing space perpendicular to the
cone axis and by integrating the lens areas 
defined by the sphere-cone intersection.
These volume integrals are rephrased with the aid
of the Byrd-Friedmann tables to Elliptic Integrals of the First, Second and Third Kind.
\end{abstract}

\maketitle

\section{On-axis. Apex inside Sphere}
\subsection{Polar coordinates}
The volume of the intersection of a cone which 
contains the sphere center on its (extended) rotation axis 
is computed in this section in tutorial fashion.
The principle parameters are the sphere radius $R$, the distance
$d$ between sphere center and cone apex, and $0\le\varphi\le \pi/2$,
half the field-of-view angle of the cone.
\begin{figure}
\includegraphics[scale=0.35]{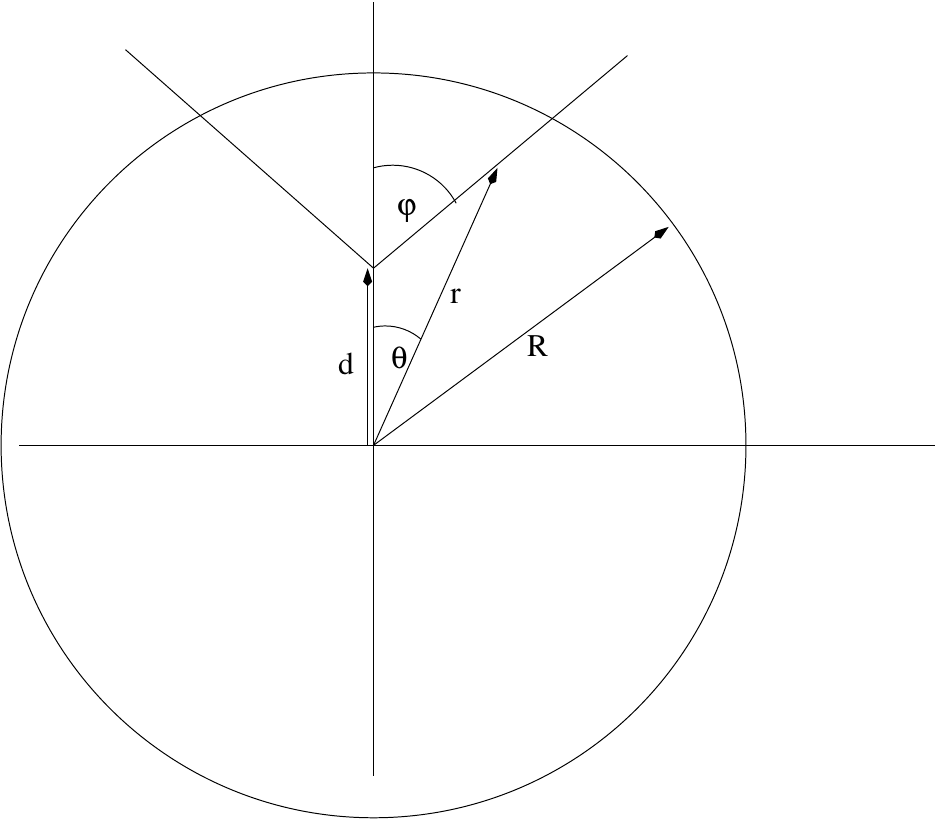}
\caption{Spherical coordinates $r$ and $\theta$ of a sphere of radius $R$
intersecting a cone with apex at a distance $d$ from the sphere center.}
\label{fig.vex}
\end{figure}
A spherical coordinate system may be defined, centered at the sphere center,
with radial coordinate $r$, azimuth angle $\phi$, polar angle $\theta$
and Jacobian $r^2\sin\theta$ \cite[4.603]{GR}.
\begin{eqnarray}
x &=& r \sin\theta \cos\phi; \\
y &=& r \sin\theta \sin\phi; \\
z &=& r \cos\theta 
\end{eqnarray}
(Some authors measure $\theta$ as a latitude from the equator and
the factor $\cos\theta$ appears in the Jacobian \cite[(1.7.9.3)]{Bronstein2}.)

\subsection{Cylinder coordinates}\label{sec.cyl}

A cylindrical coordinate system may also be used, centered at the cone apex,
with $x$ and $y$ spanning the horizontal plane,
radial coordinate $r=\sqrt{x^2+y^2}$, azimuth angle $\phi$, vertical coordinate $z$ of the cone axis,
and Jacobian $r$
\begin{eqnarray}
x &=& r \cos\phi; \\
y &=& r \sin\phi; 
\end{eqnarray}

\subsection{Volume}

\begin{figure}
\includegraphics[scale=0.35]{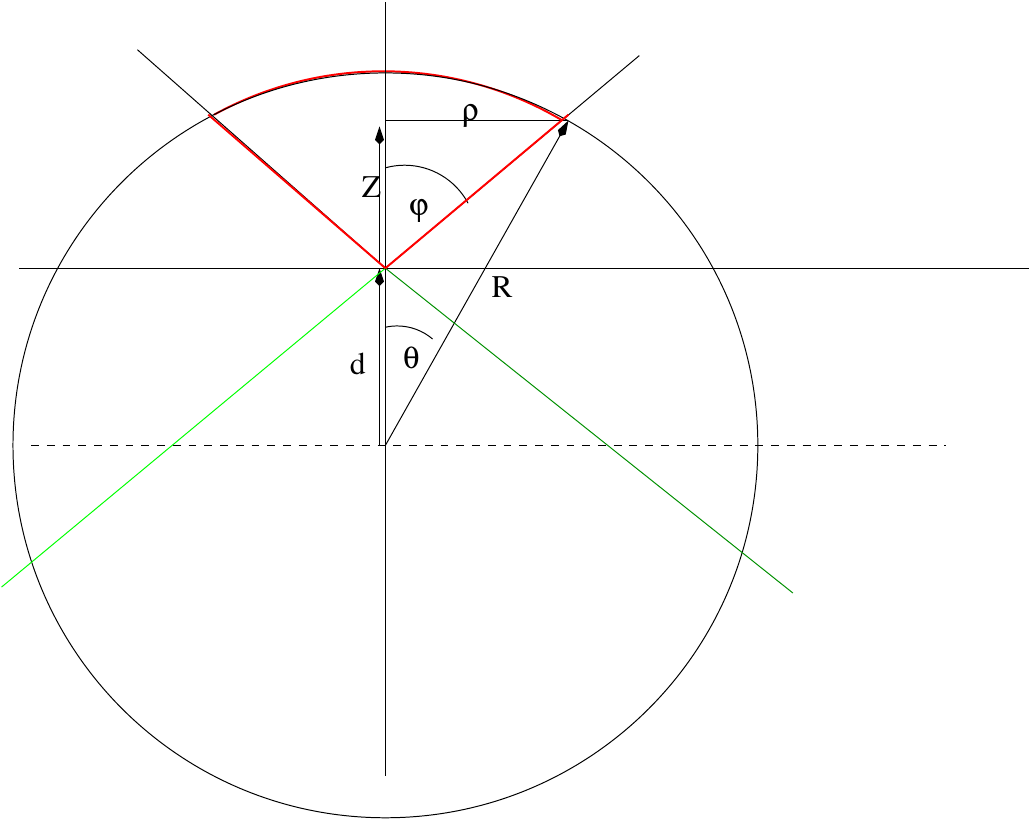}
\caption{Cylinder coordinates $z$ and $\rho$ of a sphere of radius $R$
intersecting a cone with apex at a distance $d>0$ from the sphere center. Secondary cone shell in green.}
\label{fig.vex2}
\end{figure}
\begin{figure}
\includegraphics[scale=0.4]{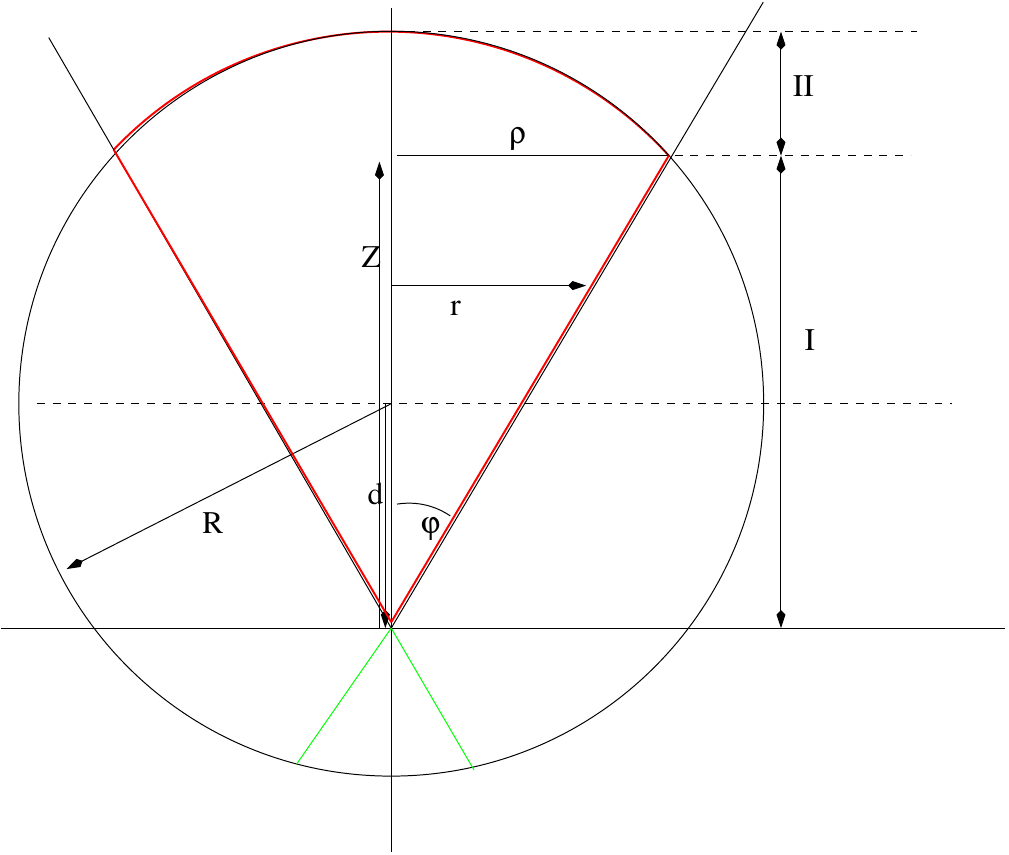}
\caption{An on-axis sphere of radius $R$
intersecting a cone with apex at a distance $d<0$ from the sphere center.}
\label{fig.incl}
\end{figure}

The volume of the intersection of the cone and sphere
is essentially the volume of the cone delimited by the horizontal plane where
the cone meets the sphere surface, plus the volume of the sphere cap for larger
values of $z$.
We define the Cartesian coordinate system such that the cone apex is at $z=0$ and
the sphere center at $z=-d$;
the sign of $d$ indicates whether the the cone apex is above or below 
the equatorial plane of the sphere.
In Figure \ref{fig.vex2} the integral starts at $z=0$ with cone radius
$r=0$ and covers the range up to $z=Z$ where the cone radius has grown linearly to $\rho$,
and continues in the range $z=Z$ up to $z=R-d$ with a radius $r$
shrinking from $\rho$ down to
zero delimited by the sphere surface.
The volume of the cone-sphere intersection is
\begin{multline}
V^{(i)}(R,d,\varphi) 
= 
\iiint dx dy dz
=\int_0^{2\pi} d\phi \int_0 r dr \int_0^{R-d} dz
=2\pi \int_0^{R-d} dz \int_0 r dr.
\label{eq.V}
\end{multline}
The transitory value $Z$ is calculated by emitting a ray from the
cone apex into the direction of $\varphi$ with coordinates $(t\sin\varphi,t\cos\varphi)$
in the cross section of Figure \ref{fig.incl}, where $t$ is the Euclidean distance
from the cone apex, and equating these with a coordinates
of a point on the sphere  with polar angle $\theta$, $(R\sin\theta,R\cos\theta-d)$.
If the $z$-coordinate of the apex is negative like in Figure \ref{fig.incl}, $d$ is smaller than zero
and the cone contains the sphere center.
Solving for $t$ gives $t=R\sin\theta/\sin\varphi$, then
\begin{equation}
R\sin(\varphi-\theta)=d\sin\varphi.
\label{eq.t0}
\end{equation}
\begin{rem}
This is also the sine equation for the
plane triangle with sides $d$, and $R$ and interior angles
$\theta$, $\pi-\varphi$ and $\varphi-\theta$ 
\cite[(3.88)]{Bronstein3}\cite[4.3.148]{AS}. 
\end{rem}
Solving for $\theta$ yields
\begin{equation}
\theta = \varphi-\arcsin(d\sin\varphi/R).
\label{eq.theta}
\end{equation}
There are two possible branches of the $\arcsin$ in this equation, the principle
value and its $\pi$-complement. The one to be taken is the one where the variable $t$
remains positive.
The other is associated with the branch of the double cone that points with its axis into the opposite direction, 
the green lines
in Figures \ref{fig.vex2} \and \ref{fig.incl}.
 The
projection on the vertical axis is
\begin{multline}
Z+d=R\cos\theta 
= R\cos[\varphi-\arcsin(d\sin\varphi/R)]
\\
= 
\cos\varphi\sqrt{R^2-(d\sin\varphi)^2}
+d\sin^2\varphi
.
\label{eq.Z}
\end{multline}

For $\varphi <\pi/2$ the intersection
of cone and sphere is at  $Z>0$; for the ``stretched'' cone where $\varphi>\pi/2$ the
intersection is at $Z<0$. The simplest strategy for the stretched cone is
to define a complementary volume which is inside the sphere but outside the cone,
\begin{equation}
\bar V^{(i)}(R,d,\varphi) \equiv \frac43 \pi R^3-V^{(i)}(R,d,\varphi) ,
\end{equation}
and to handle the cases $\varphi>\pi/2$ by flipping the cone axis upside-down:
\begin{equation}
V^{(i)}(R,d,\varphi) = \frac43 \pi R^3-V^{(i)}(R,-d,\pi-\varphi),\quad \pi/2< \varphi \le \pi.
\label{eq.Vstret}
\end{equation}
Splitting \eqref{eq.V} into the two regions $z\lessgtr Z$ yields
\begin{multline}
V^{(i)}(R,d,\varphi) 
=2\pi \int_0^{R-d} dz \int_0 r dr
\\
=2\pi \int_0^Z dz \int_0^{z\tan\varphi} r dr
+2\pi \int_Z^{R-d} dz \int_0^{\sqrt{R^2-(z+d)^2}} r dr
\\
=\pi \int_0^Z dz z^2\tan^2\varphi
+\pi \int_Z^{R-d} dz (R^2-(z+d)^2)
\\
=\frac13 \pi \tan^2\varphi Z^3
+\frac13 \pi (2R+Z+d)(R-Z-d)^2,\quad 0\le \varphi \le \pi/2.
\label{eq.vi}
\end{multline}
\begin{rem}
Equating $z\tan\varphi=\sqrt{R^2-(z+d)^2}$
in the upper limits and solving for $z$ gives again \eqref{eq.Z}.
\end{rem}
Figures \ref{fig.vex2} and \ref{fig.incl} show $\tan\varphi=\rho/Z$, so the
first term is the well-known volume 
\begin{equation}
V_\Delta(\rho,Z)=\frac13 \pi \rho^2 Z
\label{eq.cone}
\end{equation}
of the cone with
base radius $\rho$ and height $Z$ 
\cite[3.151]{Bronstein3} in region I of Figure \ref{fig.incl}, the second term is the
well-known volume 
\begin{equation}
V_\frown(R,h)=\frac13 \pi h^2(3R-h)
\label{eq.cap}
\end{equation}
of the spherical cap 
of thickness $h=R-Z-d$ \cite[3.162]{Bronstein3}, region II.

Limiting cases are 
\begin{itemize}
\item
the cone with apex at the center of the sphere with $Z=R\cos\varphi$:
\begin{equation}
V^{(i)}(R,0,\varphi)=\frac23\pi R^3(1-\cos\varphi)
\end{equation}
which is $R^3/3$ times the solid angle covered by the cone,
\item
the volume of the spherical cap 
\cite[3.162]{Bronstein3}
\begin{equation}
V^{(i)}(R,d,\pi/2)=\frac13\pi (2R+d)(R-d)^2.
\end{equation}
\item
the apex at the sphere surface with $Z=R[1+\cos(2\varphi)]$,
\begin{equation}
V^{(i)}(R,-R,\varphi)=\frac43\pi R^3 \sin^2\varphi(1+\cos^2\varphi) 
\end{equation}
such that for $\varphi\to \pi/2$ the sphere is entirely inside
the cone and the volume is the entire sphere volume $4\pi R^3/3$.
\item
the complementary spherical cap 
\begin{equation}
V^{(i)}(R,-d,\pi/2)
=
\frac43\pi R^3
-\frac13\pi (2R+d)(R-d)^2
=
\bar V^{(i)}(R,d,\pi/2)
.
\end{equation}
\end{itemize}

\section{On-axis. Apex outside the Sphere} 
\begin{figure}
\includegraphics[scale=0.4]{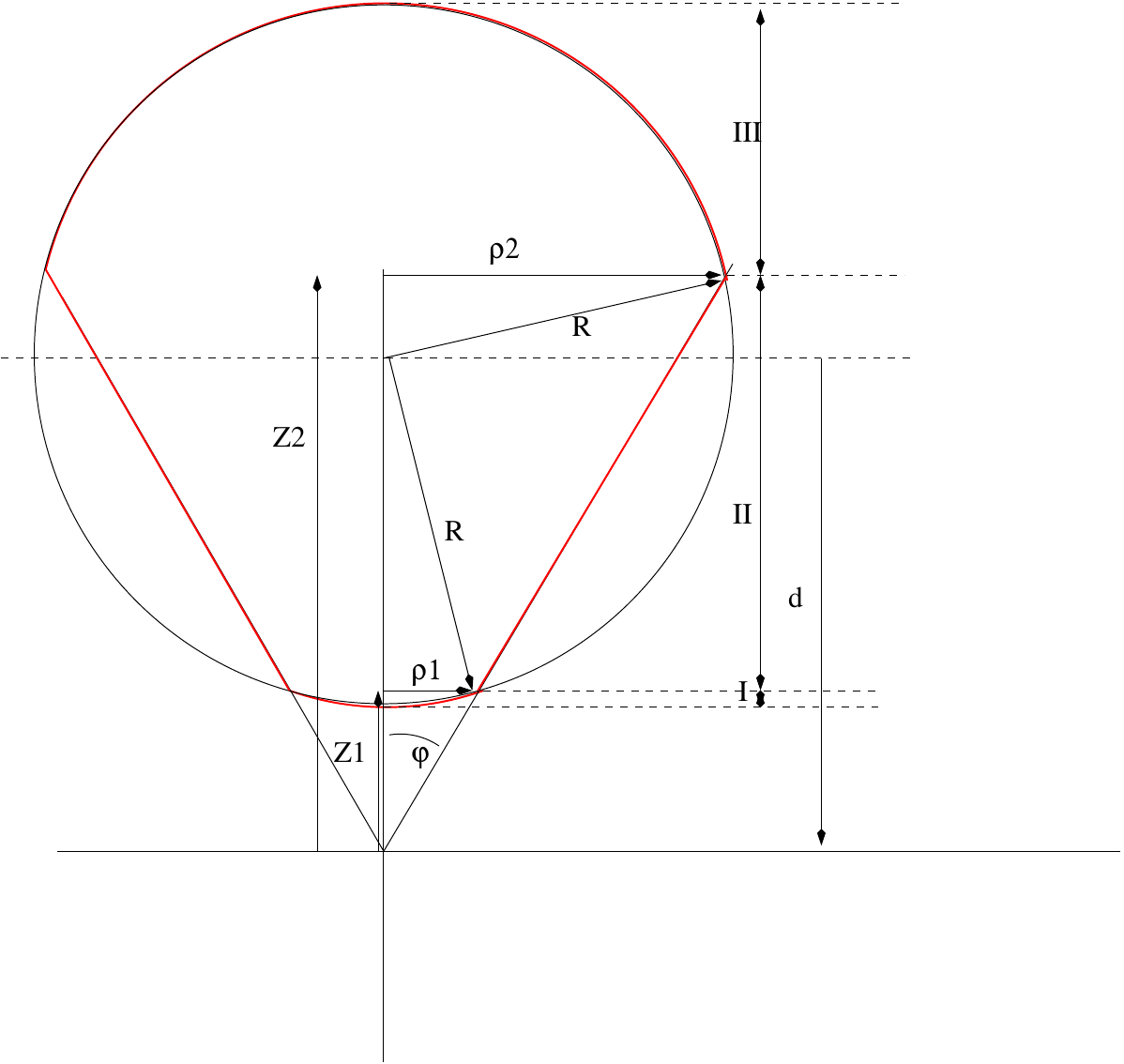}
\caption{Cone apex outside the On-axis Sphere. $d< -R$.
}
\label{fig.outs}
\end{figure}

If the cone apex is outside the sphere, $d<-R$, the cone (projection)
intersects the sphere at a near point characterized by (projected) cylinder
coordinates $Z_1,\rho_1$ and a far point $Z_2,\rho_2$ as sketched in Figure \ref{fig.outs} \cite{SchmidCGS20}.
In the figure the polar angle for the far point is $\theta< \pi/2$, but for
sufficiently small $d$ and sufficiently large $\varphi$ the far point may
actually be located below the horizontal line.

We follow the strategy to obtain the volume of the intersection as 
the sum of the volume of the sphere cap of radius $\rho_1$ at the
south pole, a truncated cone of radii $\rho_1$ and $\rho_2$, and a sphere
cap of radius $\rho_2$ at the north pole.

The ray from the apex into the direction $\varphi$ meets the
circle (projection) fixed by the condition that the distance to 
the sphere center be $R$. The same analysis that led to \eqref{eq.Z}
requires
\begin{multline}
Z_{1,2}+d
= R\cos[\varphi-\arcsin(d\sin\varphi/R)]
\\
= 
\mp
\cos\varphi\sqrt{R^2-(d\sin\varphi)^2}
+d\sin^2\varphi
.
\end{multline}
where $\arcsin$ refers to both branches of the inverse trigonometric sine,
generating two different signs of their cosines.
For sufficiently large $\varphi$, $\sin\varphi>R/|d|$, 
the cone walls may be tangential to the sphere or not intersect
the sphere at all;
then the volume of the intersection is the full $4\pi R^3/3$ of the sphere.
The corresponding (positive) radii are
\begin{equation}
\rho_{1,2}=\sqrt{R^2-(Z_{1,2}+d)^2} = Z_{1,2}\tan\varphi.
\end{equation}
The thickness of the south polar sphere cap, region I in Figure \ref{fig.outs}, is $R+Z_1+d$;
the height of the truncated cone, region II, is $Z_2-Z_1\ge 0$; 
the thickness of the north polar sphere cap, region III, is $R-Z_2-d$.
The sum of these 3 positive values is $2R$.
The sum of the volumes of the three regions is
\begin{multline}
V^{(o)}(R,d,\varphi)
=V_\frown(R,R+Z_1+d)+V_\Delta(\rho_2,Z_2)-V_\Delta(\rho_1,Z_1)+V_\frown(R,R-Z_2-d)
\\
=
\frac13 \pi \left[(R+Z_1+d)^2(2R-Z_1-d)
+(Z_2-Z_1)(\rho_1^2+\rho_1\rho_2+\rho_2^2)
+(R-Z_2-d)^2(2R+Z_2+d)
\right]
.
\label{eq.Vo}
\end{multline}
for $\sin\varphi \le R/d\le 1$, $Z_1\le Z_2$.
\begin{rem}
If $d\to -\infty$ and $\varphi\to 0$ keeping $d\sin\varphi$ constant,
the geometry approaches a sphere intersected with a circular cylinder
\cite{BoersmaPKAW64,LamarcheCPC59,SarragaCVGI22}.
\end{rem}

\section{Off-axis. Apex inside Sphere} \label{sec.offin}
\begin{figure}
\includegraphics[scale=0.4]{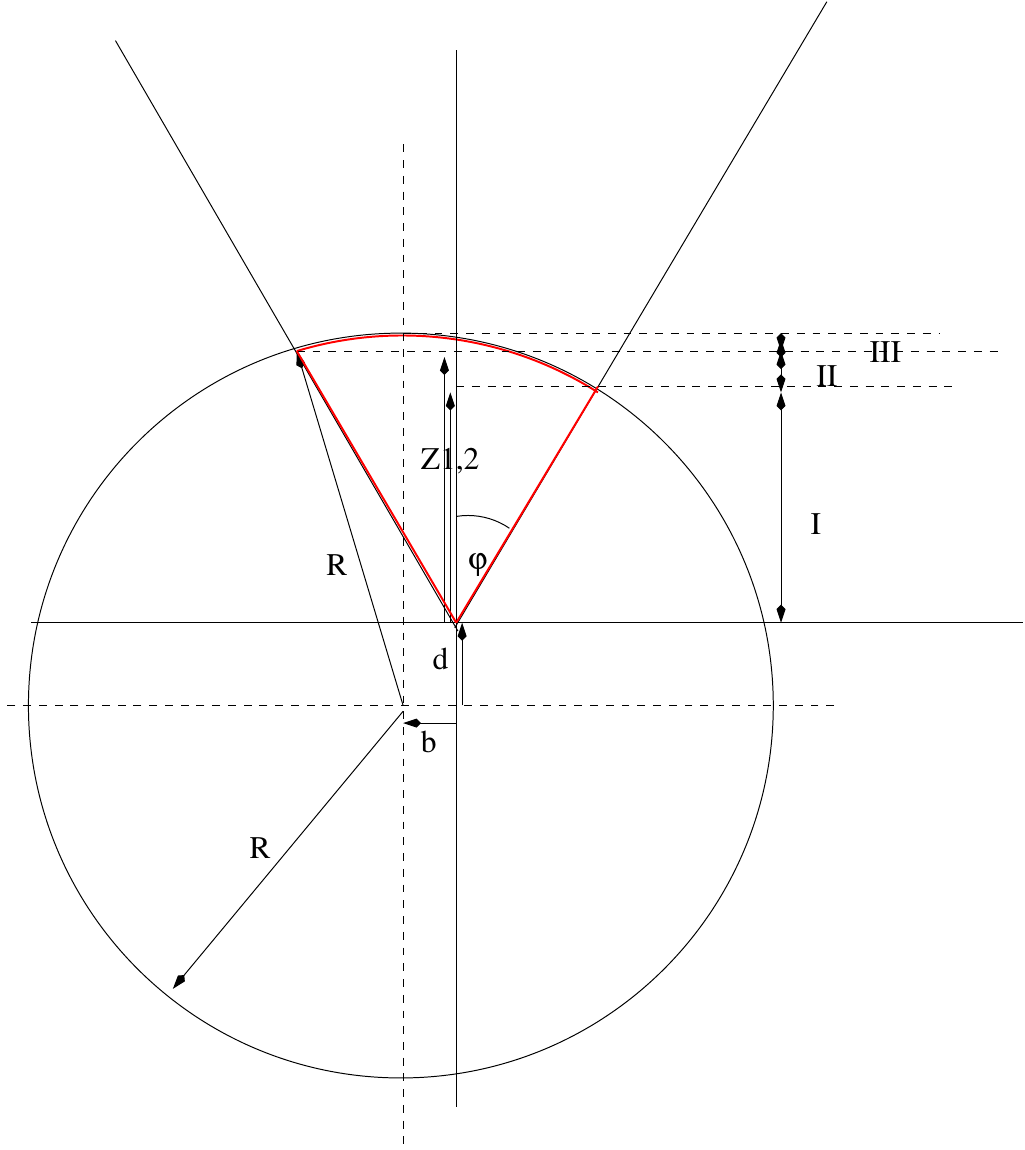}
\caption{Cone apex inside the off-axis sphere. 
Impact parameter $b$ and signed distance $d$ between apex and the sphere's equatorial plane.
}
\label{fig.incloff}
\end{figure}

The more general geometry shows a sphere center that has a nonzero distance
$b>0$, the impact parameter
of particle physics, from the cone axis as in Figure \ref{fig.incloff}.
The distance $d$ between the apex and the sphere center, measured
along the sphere axis orthogonal to $b$, may be positive or negative.
\begin{rem}
One can always find a rectangular coordinate system 
in which the cone axis and the sphere center are in the $x-z$ plane,
such that there are 4 relevant parameters: $R$, $\varphi$, $b$ and $d$.
Given a generic set of $\vec S$ (location of the sphere center), $\vec C$
(location
of the cone apex), and $\vec a$ (vector from the apex along the cone axis),
the reduction to the principal parameters $b$ and $d$ is given
by the projection of $\vec S$ on $\vec a$.
$\vec C$ is reached from $\vec S$ via $\vec S-\vec b +t\vec a=\vec C$
with $t$ measuring distances along the cone axis. 
$\vec b\perp \vec a$ yields 
in an intermediate step $t=(\vec C-\vec S)\cdot \vec a/(\vec a\cdot \vec a)$.
Then $b$ follows as $b=|\vec b|=|\vec S-\vec C+t\vec a|$,
and finally $d=t|\vec a|$.
\end{rem}

The apex is in the sphere while $d^2+b^2\le R^2$; otherwise the analysis of Section \ref{sec.offout} takes over.
We place the apex at the center of the cylinder coordinates; the
cone wall is defined by
\begin{equation}
z=r\cot\varphi; \quad z\ge 0
\label{eq.zoff}
\end{equation}
and the sphere surface by
\begin{equation}
(x+b)^2+y^2+(z+d)^2=R^2.
\end{equation}
This fixes the sign of $d$: $d$ is positive if the cone apex is \emph{above} the
altitude of the sphere center.
Section \ref{sec.cyl} turns this into
\begin{equation}
r^2+b^2+2rb\cos\phi+(z+d)^2=R^2
\label{eq.offCyc}
.
\end{equation}
In the projection, Figure \ref{fig.incloff},
a lower value $Z_1$ and a higher value $Z_2$ of the latitudes appear where
sphere and cone intersect.
\begin{rem}
\label{rem.critr}
On the curve of the intersection one can eliminate $z$ from
\eqref{eq.offCyc} inserting \eqref{eq.zoff}:
\begin{equation}
r^2+b^2+2rb\cos\phi+(r\cot\varphi+d)^2=R^2.
\label{eq.offCyc2}
\end{equation}
This is the projection of the intersection curve on the horizontal plane (perpendicular to the cone
axis).
The first derivative of \eqref{eq.offCyc2}
along the curve of $r(\phi)$, with $r'\equiv \partial r/\partial \phi$, is
\begin{equation}
2rr'+2r'b\cos\phi-2rb\sin\phi+2(r\cot\varphi+d)r'=0.
\end{equation}
The maximum and minimum distance of this planar curve from the
cone axis are set by $r'=0$ which leads with this equation
to $\sin\phi=0$; so it is correct to
locate the points of $Z_{1,2}$ in the cross section of the $x-z$-plane.
\end{rem}

\begin{rem}
We only consider angles $\varphi\le \pi/2$ here. The cases of overstretched 
cones are forwarded to \eqref{eq.Vstret} with flipped sign of $d$.
\end{rem}

The locations of $Z_{1,2}$ are calculated in spherical coordinates centered at the cone apex
sending a ray of parameter
$t$ from the apex into the direction $\varphi$ along $(x,z)$ coordinates
$(t\sin\varphi,t\cos\varphi)$. Another ray from the sphere center
to the sphere surface has coordinates $(R\sin\theta-b,R\cos\theta-d)$
where $\theta$ is the angle measured from the center of the sphere.
They meet at
\begin{equation}
t\sin\varphi = R\sin\theta-b;\quad
t\cos\varphi = R\cos\theta-d.
\end{equation}
Introducing unitless
\begin{equation}
\hat d\equiv d/R;\quad \hat b\equiv b/R
\end{equation}
and eliminating $t$ gives
\begin{equation}
\sin(\varphi -\theta)
= \hat d\sin\varphi -\hat b\cos\varphi,
\end{equation}
\begin{equation}
\theta=
\varphi
-\arcsin(\hat d\sin\varphi -\hat b\cos\varphi)
\label{eq.thetab}
\end{equation}
consistent with the on-axis equation \eqref{eq.t0} as $b\to 0$.
\begin{rem}
There are two branches of the $\arcsin$, selecting the cone and
its twin at the same apex but the axis pointing into the opposite direction. 
\end{rem}
Back to the cylinder coordinates centered at the apex the projections are
\begin{equation}
Z_{1,2}+d=R\cos\theta
\end{equation}
switching to $-\varphi$ in \eqref{eq.thetab} to obtain the second value:
\begin{eqnarray}
Z_1+d&=&R\cos[\varphi-\arcsin(\hat d\sin\varphi-\hat b \cos\varphi];\label{eq.Z1}\\
Z_2+d&=&R\cos[-\varphi-\arcsin(-\hat d\sin\varphi-\hat b \cos\varphi].
\end{eqnarray}

The limits of the radii of the cylinder coordinates are
\begin{equation}
\rho_{1,2}=Z_{1,2}\tan\varphi.
\label{eq.rhoofZ}
\end{equation}

\begin{rem}
An alternative to this vector-algebra
is to continue from Remark \ref{rem.critr} with an algebraic approach:
Since the critical values are where $\phi=0,\pi$, one may insert
$\cos\phi=\pm 1$ into \eqref{eq.offCyc2}, solve the two quadratic equations for $r$
in the range $r>0$, call them
$\rho_{1,2}$ and insert these back into \eqref{eq.zoff} to get $Z_{1,2}$.
\end{rem}
\begin{rem}
In the limit $b=0$ both solutions of \eqref{eq.offCyc2} are the same because
the term with $\cos\phi$ vanishes from the equation.
\end{rem}

The intersecting volume $V^{(i)}(R,d,b,\varphi)=V_I+V_{II}+V_{III}$
contains three subvolumes if the $z$-axis is chopped by horizontal
planes into slices. 
\begin{enumerate}
\item
For $0<z< Z_1$, region I in Figure \ref{fig.incloff}, it
is the volume of the cone of height $Z_1$ and base radius $\rho_1\equiv Z_1\tan\varphi$;
this subvolume has the value $V_I=V_\Delta(\rho_1,Z_1)=\frac13 \pi \rho_1^2 Z_1$.
\item
 For $Z_1<z<Z_2$, region II, the limits of $r$ are partially set by the cone wall
and partially by the sphere surface.
\item
For $Z_2<z<R-d$ 
the limits of $r$, region III, are set by the sphere surface as a function of 
the azimuth $\phi$.
It is a sphere cap of base
radius $\rho_2-b$ and thickness $h=R-d-Z_2$ thinking
in terms of polar coordinates at the sphere center.
This portion of the volume of the North Polar Cap is
$V_{III}=V_\frown(R,h)$.
For sufficiently small $\varphi$, $Z_2\tan\varphi < b$, the sphere cap is outside
the cone
and there is no contribution of region III to the volume. Then the left branch of the red
line of the cone in Figure \ref{fig.excloff} hits the circle 
before the red arc on the circle reaches
the north pole.
[This type of argument indicating the presence of sphere caps in the 
volume of intersection will recur many times in Section \ref{sec.offout}.]
\end{enumerate}

This leaves to find a formula for the volume $V_{II}$ of the intermediate region (II).
For constant $z$ the interior of the cone defines a circle
of radius $z\tan\varphi$ at the center, and the interior of the sphere
defines a circle of radius $\sqrt{R^2-(z+d)^2}$ displaced by $b$.
In region (2) none of both is entirely inside the other, 
and the horizontal slice defines an overlap area $A$ as discussed
in Appendix \ref{app.2c}.

The volume integral gathers the cross section areas \eqref{eq.A} along the $z$-direction,
\begin{equation}
V_{II}=\int_{Z_1}^{Z_2} A\left(z\tan\varphi, \sqrt{R^2-(z+d)^2},b\right)dz,
\label{eq.vcosph}
\end{equation}
which splits into three terms overlapping two sectors and subtracting one triangle,
\begin{equation}
V_{II}=v_1 +v_2-v_\Delta.
\label{eq.v3incl}
\end{equation}
The ordinate sections of \eqref{eq.lensOrd} are for these $z$-dependent radii
\begin{equation}
x_1=\frac{1}{2b}(R^2-\frac{z^2}{\cos^2\varphi}-2zd-d^2-b^2),
\end{equation}
\begin{equation}
x_2=\frac{1}{2b}(R^2-\frac{z^2}{\cos^2\varphi}-2zd-d^2+b^2).
\end{equation}
To consolidate the algebraic representation it is useful to introduce the
dimensionless variable
\begin{equation}
\hat z\equiv z/(R\cos\varphi).
\label{eq.zhat}
\end{equation}

\subsection{Cone Sector} 
The contribution $v_1$ to \eqref{eq.v3incl} along the vertical direction induced by the area of \eqref{eq.A} is 
\begin{multline}
v_1
=\int dz \alpha_1 r_1^2
=\int dz 
r_1^2
\arccos\frac{-x_1}{r_1} 
\\
=\int dz 
z^2\tan^2\varphi 
\arccos\frac{-\frac{1}{2b}(R^2-\frac{z^2}{\cos^2\varphi}-2zd-d^2-b^2)}{z\tan\varphi } 
\\
=R^3\cos\varphi \sin^2\varphi
\int d\hat z \hat z^2 \arccos\frac{-1+\hat z^2+2\cos\varphi \hat d\hat z+\hat d^2+\hat b^2}{2\hat b\hat z\sin\varphi } 
.
\end{multline}
With partial integration and $\frac{d}{dx} \arccos x= -1/\sqrt{1-x^2}$ for the outer derivative:
\begin{multline}
\int d\hat z \hat z^2 \arccos\frac{-1+\hat z^2+2 \hat d\hat z+\hat d^2+\hat b^2}{2\hat b\hat z\sin\varphi }
\\
=
\frac13 \hat z^3 \arccos[\ldots]-\frac13 \int d \hat z \hat z^3 
\\
\frac{-\hat z^2-1+\hat d^2+\hat b^2}{
\hat z^2\hat b\sin\varphi
\sqrt{
-\frac{
(2\cos\varphi \hat d \hat z-1+\hat z^2+\hat b^2+\hat d^2+2\hat b\sin\varphi \hat z)
(2\cos\varphi \hat d \hat z-1+\hat z^2+\hat b^2+\hat d^2-2\hat b\sin\varphi \hat z)
}{\hat b^2\sin^2\varphi \hat z^2}
}
}
\\
=
\frac13 \hat z^3 \arccos[\ldots]+\frac13 \int d \hat z \hat z^3 
\\
\frac{\hat z^2+1-\hat d^2-\hat b^2}{
\hat z
\sqrt{
-
(2\cos\varphi \hat d \hat z-1+\hat z^2+\hat b^2+\hat d^2+2\hat b\sin\varphi \hat z)
(2\cos\varphi \hat d \hat z-1+\hat z^2+\hat b^2+\hat d^2-2\hat b\sin\varphi \hat z)
}
}
\\
=
\frac13 \hat z^3 \arccos[\ldots]+\frac13 \int d \hat z 
\frac{\hat z^2(\hat z^2+1-\hat d^2-\hat b^2)}{
\sqrt{
-
(\hat z -\hat z_1^+)
(\hat z -\hat z_1^-)
(\hat z -\hat z_2^+)
(\hat z -\hat z_2^-)
}
}
,
\label{eq.ellcone}
\end{multline}
where the roots of the two quadratic polynomials in the denominator are
\begin{equation}
\hat z_{1}^{\pm} = 
-(\hat d \cos\varphi +\hat b \sin\varphi)
\pm \sqrt{(\hat d \cos\varphi +\hat b \sin\varphi)^2+1-\hat d^2-\hat b^2} ; \, \hat z_1^{-}<0<\hat z_1^{+}
\label{eq.z1}
\end{equation}
\begin{equation}
\hat z_2^{\pm} = 
-(\hat d \cos\varphi -\hat b \sin\varphi)
\pm \sqrt{(\hat d \cos\varphi -\hat b \sin\varphi)^2+1-\hat d^2-\hat b^2} ;
\, \hat z_2^{-}<0<\hat z_2^{+}
\label{eq.z2}
\end{equation}
$\hat z_1^+$ is the scaled version of $Z_1$, $\hat z_2^+$ is the scaled version of $Z_2$.
If $\hat d^2+\hat b^2<1$ the apex is inside the sphere of radius $R$.
The integral is reduced as detailed in Appendix \ref{app.Csec}.

\subsection{Sphere Sector} 
The contribution $v_2$ to \eqref{eq.v3incl} along the vertical direction induced by the area of \eqref{eq.A} is in
the scaled variable \eqref{eq.zhat}
\begin{multline}
v_2
=\int dz \alpha_2 r_2^2
=\int dz r_2^2 \arccos\frac{x_2}{r_2} 
\\
=\int dz 
[R^2-(z+d)^2] 
\arccos\frac{\frac{1}{2b}(R^2-\frac{z^2}{\cos^2\varphi}-2zd-d^2+b^2)}{\sqrt{R^2-(z+d)^2}} 
\\
=R^3\cos\varphi 
\int d\hat z 
[1-(\hat z\cos\varphi +\hat d)^2]
\arccos\frac{1-\hat z^2-2\cos\varphi \hat d\hat z-\hat d^2+\hat b^2}{2\hat b\sqrt{1-(\cos\varphi \hat z+\hat d)^2}} 
.
\end{multline}
By partial integration
\begin{multline}
\cos\varphi\int d\hat z 
\left[1-(\hat z\cos\varphi +\hat d)^2\right]
\arccos\frac{1-\hat z^2-2\cos\varphi \hat d\hat z-\hat d^2+\hat b^2}{2\hat b\sqrt{1-(\cos\varphi \hat z+\hat d)^2}} 
\\
=
\cos\varphi \left[\hat z-\frac{(\hat z\cos\varphi+\hat d)^3}{3\cos\varphi}\right]
\arccos[\ldots]
+
\cos\varphi \int d\hat z 
\left[\hat z-\frac{(\hat z\cos\varphi+\hat d)^3}{3\cos\varphi}\right]
\\
\times \frac{1}{\sqrt{1-\frac{(1-\hat z^2-2\cos\varphi \hat d\hat z-\hat d^2+\hat b^2)^2}{4\hat b^2 [1-(\cos\varphi \hat z+\hat d)^2]}
}}
\\
\times \frac{1}{2\hat b\sqrt{1-(\cos\varphi \hat z+\hat d)^2}}\left[
-2\hat z-2\hat d\cos\varphi
+\frac{(1-\hat z^2-2\hat d\cos\varphi \hat z-\hat d^2+\hat b^2)(\cos\varphi \hat z+\hat d)\cos\varphi}{1-(\cos\varphi \hat z+\hat d)^2}
\right]
\\
=
\cos\varphi \left[\hat z-\frac{(\hat z\cos\varphi+\hat d)^3}{3\cos\varphi}\right]
\arccos[\ldots]
+
\cos\varphi \int d\hat z 
\left[\hat z-\frac{(\hat z\cos\varphi+\hat d)^3}{3\cos\varphi}\right]
\\
\times \frac{1}{\sqrt{4\hat b^2\left[1-(\cos\varphi \hat z+\hat d)^2\right]-(1-\hat z^2-2\cos\varphi \hat d\hat z-\hat d^2+\hat b^2)^2}
}
\\
\times
\left[
-2\hat z-2\hat d\cos\varphi
+\frac{(1-\hat z^2-2\hat d\cos\varphi \hat z-\hat d^2+\hat b^2)(\cos\varphi \hat z+\hat d)\cos\varphi}{1-(\cos\varphi \hat z+\hat d)^2}
\right]
\\
=
\left[\hat z\cos\varphi -\frac{(\hat z\cos\varphi+\hat d)^3}{3}\right]
\arccos[\ldots]
+
\int d\hat z 
\left[\hat z\cos\varphi -\frac{(\hat z\cos\varphi+\hat d)^3}{3}\right]
\\
\times
\left[
-2(\hat z+\hat d\cos\varphi)
+\cos\varphi\frac{(1-\hat z^2-2\hat d\cos\varphi \hat z-\hat d^2+\hat b^2)(\hat z\cos\varphi+\hat d)}
{
(1+\cos\varphi \hat z+\hat d)
(1-\cos\varphi \hat z-\hat d)
}
\right]
\\
\times \frac{1}{\sqrt{
(\hat z-\hat z_1^{+})
(\hat z-\hat z_1^{-})
(\hat z_2^+-\hat z)
(\hat z-\hat z_2^{-})
}
}
\label{eq.ellsph}
\end{multline}
with $\hat z_{1,2}^{\pm}$ given by \eqref{eq.z1}--\eqref{eq.z2}.
\begin{rem}
If the cone apex is on the sphere surface, $\hat d^2+\hat b^2=1$ and two of the four polynomial roots $z_{1,2}^\pm$ 
are zero such that
 a spurious singularity $\propto 1/z$ arises from the denominator product $\surd (\hat z-\hat z_1^{+})
(\hat z-\hat z_1^{-})
(\hat z_2^+-\hat z)
(\hat z-\hat z_2^{-})$. In this case this singularity can be lifted by writing the kernel
of the integral of
the previous equation (two brackets and square root) as
\begin{multline}
\int d\hat z 
[
\frac13 \hat d(-3+2\cos^2\varphi +2\hat d^2)
+\frac13 \cos\varphi( -4+2\cos^2\varphi+5\hat d^2)\hat z
+\frac43 \cos^2\varphi  \hat d^2\hat z^2
+\frac13 \cos^3\varphi  hat z^3
\\
+\frac{\hat d^2/2+5\hat d/6-\cos^2\varphi \hat d+1/3 -2\cos^2\varphi/3}
{1+\cos\varphi \hat z +\hat d}
+\frac{-\hat d^2/2+5\hat d/6-\cos^2\varphi \hat d+1/3 -2\cos^2\varphi/3}
{1-\cos\varphi \hat z -\hat d}
]
\\
\times
\frac{1}{\sqrt{
(\hat z-\hat z_1^{-})
(\hat z_2^{-}-\hat z)
}}.
\end{multline}
\end{rem}

Then \eqref{eq.ellsph} contains first a sum of integrals of powers of $\hat z$ divided by the quartic square root
manageable by Appendix \ref{app.Csec}. In addition there is an integral
with a quadratic polynomial in $\hat z$ in the denominator, which is split
into two terms with a linear polynomial by decomposition into partial fractions.
To reduce complexity introduce an intermediate
$\bar z \equiv \hat z \cos\varphi+\hat d$, so

\begin{multline}
\cos^2 \varphi
\left[\hat z-\frac{(\hat z\cos\varphi+\hat d)^3}{3\cos\varphi}\right]
\\
\times
\left[
-2(\hat z+\hat d \cos\varphi)+\cos\varphi
\frac{(1-\hat z^2-2\hat d\cos\varphi \hat z-\hat d^2+\hat b^2)(\hat z\cos\varphi+\hat d)}
{
(1+\cos\varphi \hat z+\hat d)
(1-\cos\varphi \hat z-\hat d)
}
\right]
\\
=
\frac23+\hat d^2
+\cos^2\varphi[
-\frac23-\frac23 \hat b^2-\hat d^2+\frac13 \hat d ^4+\frac13 \hat b^2\hat d^2
]
\\
+\frac13 \hat d \cos\varphi \left[-3+2\hat b^2\cos^2\varphi+2\hat d^2+2\hat d^2\cos^2\varphi\right]\hat z
\\
+\frac13\cos^2\varphi \left[
-4+\cos^2\varphi+\hat b^2\cos^2\varphi
+5\hat d^2+\hat d^2\cos^2\varphi)
\right]\hat z^2
\\
+\frac43 \hat d \cos^3\varphi \hat z^3
+\frac13 \cos^4\varphi \hat z^4
\\
+\frac{(-\frac43 \hat d^2+\frac76 \hat d+\frac12 \hat d^3-\frac13)\sin^2\varphi
+(-\frac12 \hat d\hat b^2+\frac13 \hat b^2)\cos^2\varphi
 }{1-\hat z \cos\varphi-\hat d}
\\
+\frac{(-\frac43 \hat d^2-\frac13-\frac76 \hat d-\frac12 \hat d^3) \sin^2\varphi
+(\frac13 \hat b ^2+\frac12 \hat d \hat b^2)\cos^2\varphi
  }{1+\hat z\cos\varphi +\hat d}
.
\label{eq.d}
\end{multline}
This type of integrals is covered by Appendix \ref{app.Ssec}.

\subsection{Triangle} 
The radius of the lens of the intersecting planar circles is
according to \eqref{eq.rho}
\begin{equation}
\rho
=\frac{1}{2b}
\sqrt{-(R^2-d^2-b^2-2zd-z^2/\cos^2\varphi -2bz\tan\varphi)
(R^2-d^2-b^2-2zd-z^2/\cos^2\varphi +2bz\tan\varphi)}.
\end{equation}
Then the third contribution to \eqref{eq.v3incl} can be written as
\begin{multline}
v_\Delta = \int \rho b dz 
=\frac{1}{2}R^3\cos\varphi 
\int d\hat z  \\
\times 
\sqrt{-(-1+\hat d^2+\hat b^2+2\hat d\cos\varphi \hat z +\hat z^2 +2\hat b \hat z\sin\varphi)
(-1+\hat d^2+\hat b^2+2\hat d \cos\varphi \hat z+\hat z^2 -2\hat b\hat z\sin\varphi)}
\end{multline}
The square root contains a product of two quadratic polynomials in $\hat z$.
Their roots are provided by \eqref{eq.z1}--\eqref{eq.z2}.
Since $\hat z$ will be integrated in these limits, the writeup as
an Elliptic Integral is
\begin{equation}
v_\Delta
=\frac{1}{2}R^3\cos\varphi 
\int_{z_1^+}^{z_2^+} d\hat z 
\sqrt{(\hat z-\hat z_1^+)(\hat z-\hat z_1^-)
(\hat z_2^+-\hat z)(\hat z-\hat z_2^-)}.
\end{equation}
Details of the evaluation are posted in Appendix \ref{app.Atri}.

\section{Off-axis. Apex outside Sphere} \label{sec.offout}
Determining the curve of the intersection is of interest
for satellite imagery where the sphere represents the Earth and the
cone apex an orbiting satellite \cite{RuffJAM10}.

Miller \cite{MillerACMTG6} considers
five major cases of intersections which can be registered
by the number of real-valued quantities in Equations \eqref{eq.z1}--\eqref{eq.z2}. 

If all 4 values $\hat z_{1,2}^\pm$ are imaginary, the bodies either do not intersect
at all or the sphere is entirely inside the cone with intersecting volume 
of $0$ or $4\pi R^3/3$, respectively.

\begin{figure}
\includegraphics[scale=0.4]{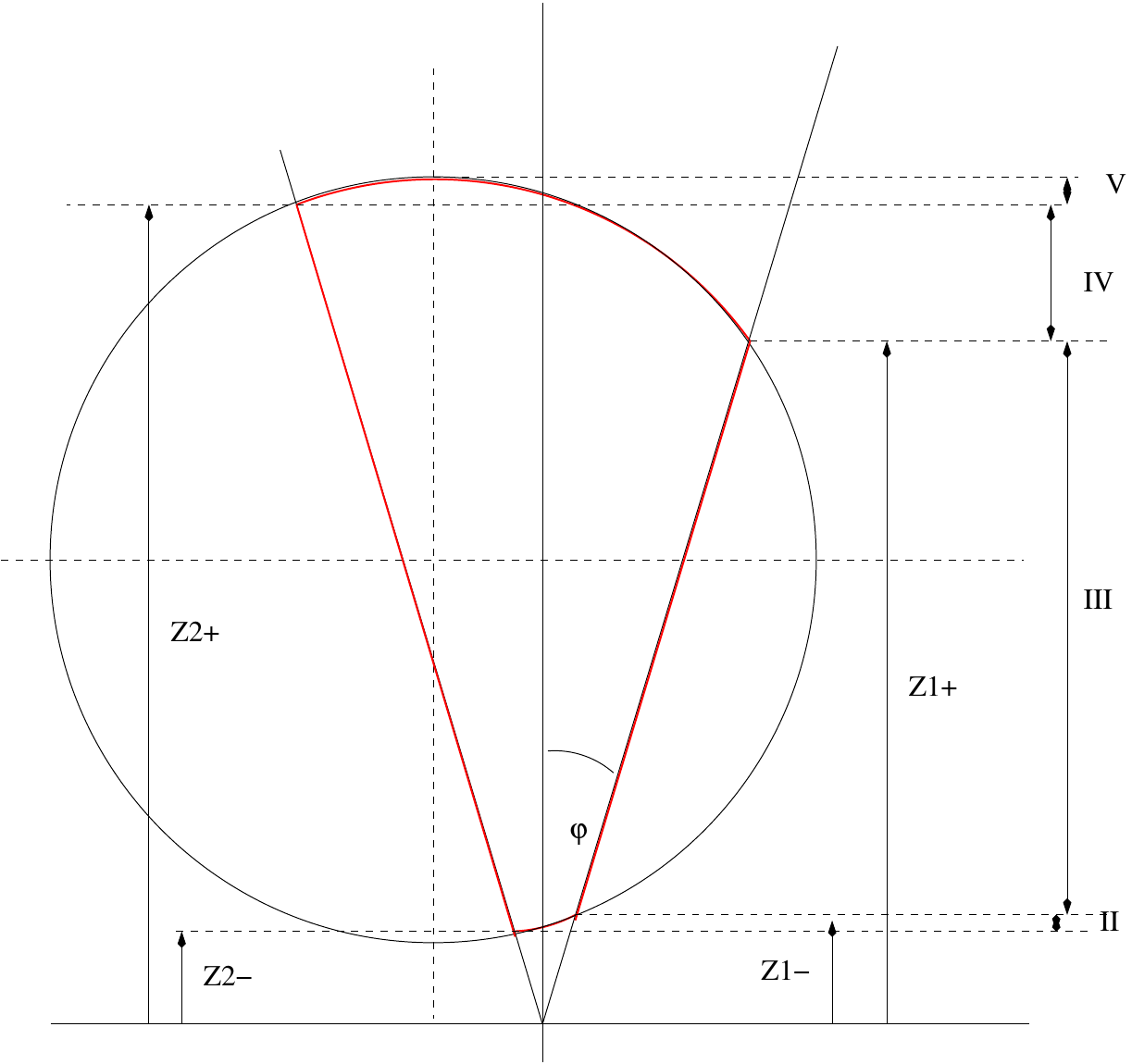}
\caption{Cone apex outside the off-axis sphere. Two-branch intersection.
Real-valued $Z_{1,2}^\pm$ and real-valued reduced $\hat z_{1,2}^\pm$.
No south pole cap in intersection.
}
\label{fig.excloff}
\end{figure}

\begin{figure}
\includegraphics[scale=0.4]{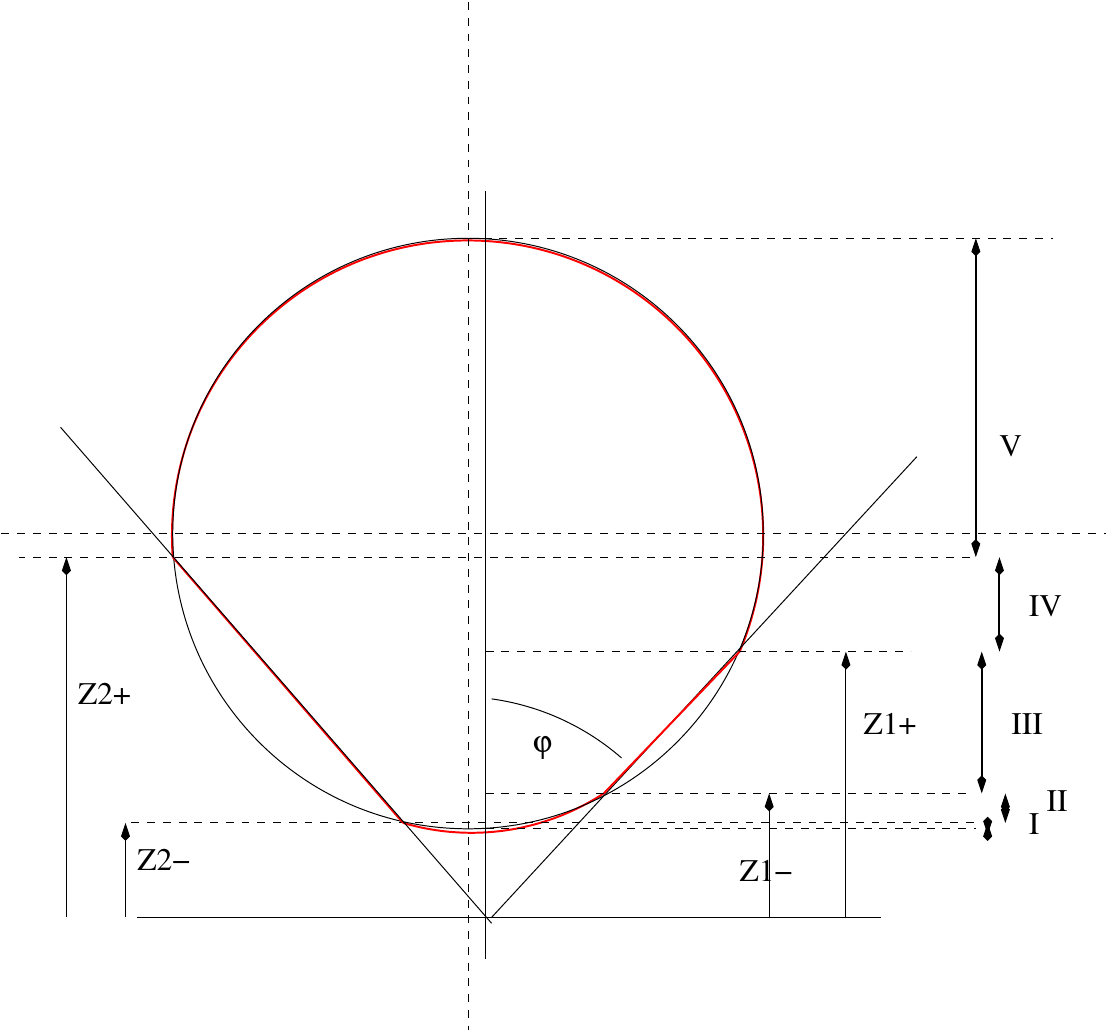}
\caption{Cone apex outside the off-axis sphere. Two-branch intersection.
Real-valued $Z_{1,2}^\pm$ and real-valued reduced $\hat z_{1,2}^\pm$.
With south pole cap in intersection.
}
\label{fig.excloffb}
\end{figure}

\begin{figure}
\includegraphics[scale=0.4]{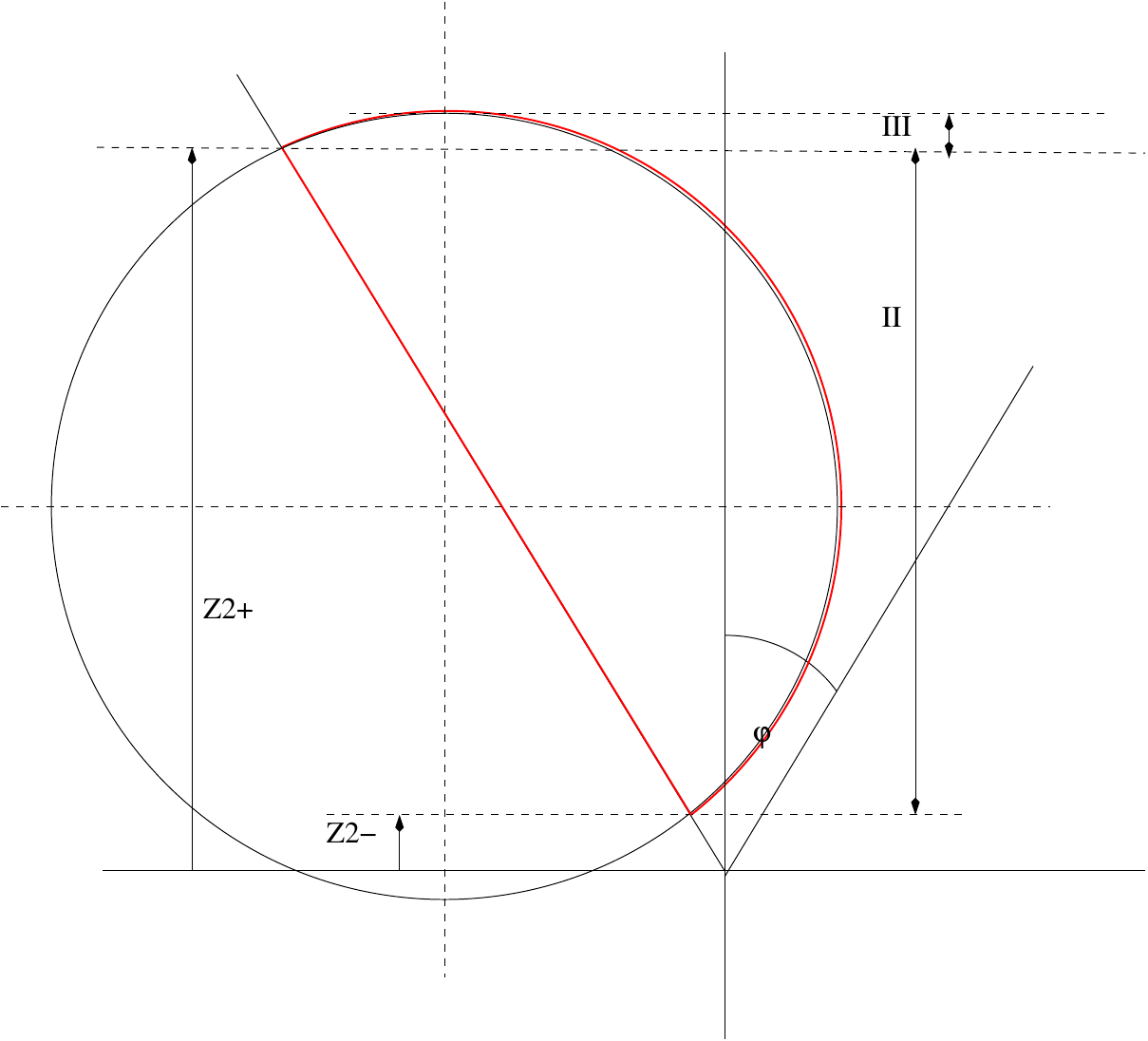}
\caption{Cone apex outside the off-axis sphere. One-branch intersection.
Real-valued $Z_2^\pm$ and real-valued reduced $\hat z_2^\pm$.
}
\label{fig.excloff2}
\end{figure}

There are essentially the two-branched case of Figure \ref{fig.excloff}
of \ref{fig.excloffb}, where all
4 values $z_{1,2}^{\pm}$ are real and the lines of intersection are two separate
quadrics, and the one-branched case of Figure \ref{fig.excloff2} where
the 2 values $z_1^{\pm}$ are imaginary and the lines of intersection is a single
quadric. 

The cases with a single tangent point where $z_1^+=z_1^-$ are just limiting
values (the Viviani case, so to speak) and not of special importance to the 
computation of the volume of the intersection of the two bodies.
\begin{rem}
$z_1^{\pm}$ are imaginary if the argument 
$(\hat d \cos\varphi +\hat b \sin\varphi)^2+1-\hat d^2-\hat b^2$=
$1-(\hat d \sin\varphi-\hat b \cos\varphi)^2$ of the square root \eqref{eq.z1} is negative.
This is the same criterion which keeps the value of the
arcsin argument in \eqref{eq.Z1} outside the interval $[-1,1]$.
\end{rem}

\subsection{Two Branches} 
In Figures \ref{fig.excloff} and \ref{fig.excloffb}, 
the intersection $V=V_I+V_{II}+V_{III}+V_{IV}+V_V$ comprises
\begin{enumerate}
\item
for sufficiently large $\varphi$ like in Figure \ref{fig.excloffb}
but not in Figure \ref{fig.excloff}, $Z_2^-\tan\varphi>b$, a sphere
cap of thickness $h=R+d+Z_2^-$ at the south pole, with volume
$
V_I = V_\frown(R,h).
$
\item
\begin{equation}
V_{II}=\int_{Z_2^-}^{Z_1^-} A\left(z\tan\varphi, \sqrt{R^2-(z+d)^2},b\right)dz
\end{equation}
\item
the truncated cone
\begin{equation*}
V_{III}=V_\Delta(\rho_1^+,Z_1^+)-V_\Delta(\rho_1^-,Z_1^-),
\end{equation*}
where $\rho$ and $Z$ are correlated by \eqref{eq.rhoofZ};
\item
\begin{equation}
V_{IV}=\int_{Z_1^+}^{Z_2^+} A\left(z\tan\varphi, \sqrt{R^2-(z+d)^2},b\right)dz
\end{equation}
\item
for sufficiently large $\varphi$, $Z_2^+\tan\varphi>b$,
a sphere cap of thickness $h=R-d-Z_2^+$ at the north pole
covering the range $Z_2^+\le z \le R-d$
with volume
$
V_{V}=V_\frown(R,h).
$
For $Z_2^+\tan\varphi<b$
the intersection at altitude $Z_2^+$
may be at the horizontal cylinder coordinate smaller than $b$;
then this sphere cap is outside the cone and does not contribute to the volume.
\end{enumerate}
The algebra for regions II and IV is the same as the algebra of computing
\eqref{eq.vcosph}. The only difference is that the numerical
order of the four real-valued reduced $\hat z$ values in 
$\sqrt{-(\hat z-z_1^+)(\hat z-z_1^-)(\hat z-z_2^+)(\hat z-z_2^-)}$
may differ such that the parameters $a$--$d$ in the appendices
need to be permuted to yield applicable Byrd-Friedman-integrals.

\subsection{One Branch} 
If only $\hat z_2^\pm$ are real-valued the volume of the
intersection $V=V_I+V_{II}+V_{III}$ contains 
\begin{enumerate}
\item
for sufficiently large $\varphi$, $Z_2^-\tan\varphi >b$ a sphere cap
of thickness $h=R+d+Z_2^-$ at the south pole, volume $V_I=V_\frown(R,h)$.
In Figure \ref{fig.excloff2} that south pole cap is not contributing.
\item
a subvolume 
\begin{equation*}
V_{II}=\int_{Z_2^-}^{Z_2^+} A\left(z\tan\varphi, \sqrt{R^2-(z+d)^2},b\right)dz
\end{equation*}
as in \eqref{eq.vcosph}.
For the main integral the analysis of Section \ref{sec.offin} remains valid, but here
$\hat z_1^+$ and $\hat z_1^-$ are two conjugated   complex values, so
in the appendices the alternative integrals with complex-conjugated roots
of the quartic polynomial are activated.
\item
for sufficiently large $\varphi$, $Z_1^+\tan\varphi >b$ a sphere cap
of thickness $h=R-d-Z_2^+$ for $Z_2^+<z<R-d$ up to the north pole,
volume $V_{III}=V_\frown(R,h)$.
\end{enumerate}
In overview, the one-branch cases have the same criteria and formulas
to include the polar caps as the two-branch cases, they have no
contribution from a truncated cone in intermediate $Z$-regions, 
and the two integrals
that depend on the limits $Z_1^{\pm}$ in the two-branch cases are
glued into a single integral covering the entire interval $[Z_2^-,Z_2^+]$. 

\section{Summary}
The volume of the intersection of cone and on-axis sphere is
given 
by \eqref{eq.vi} if the apex is inside the sphere,
and by \eqref{eq.Vo} 
and the apex is outside the sphere: sums of (truncated) cones
and sphere caps.

Answering a question of Shah \cite{ShahSIAMR30}, the
volumes with off-axis spheres have been reduced to Elliptic Integrals in
Sections \ref{sec.offin} and \ref{sec.offout} for apexes in- and outside the sphere.

\newpage 
\appendix
\section{Two Intersecting Circles} \label{app.2c}
\begin{figure}
\includegraphics[scale=0.45]{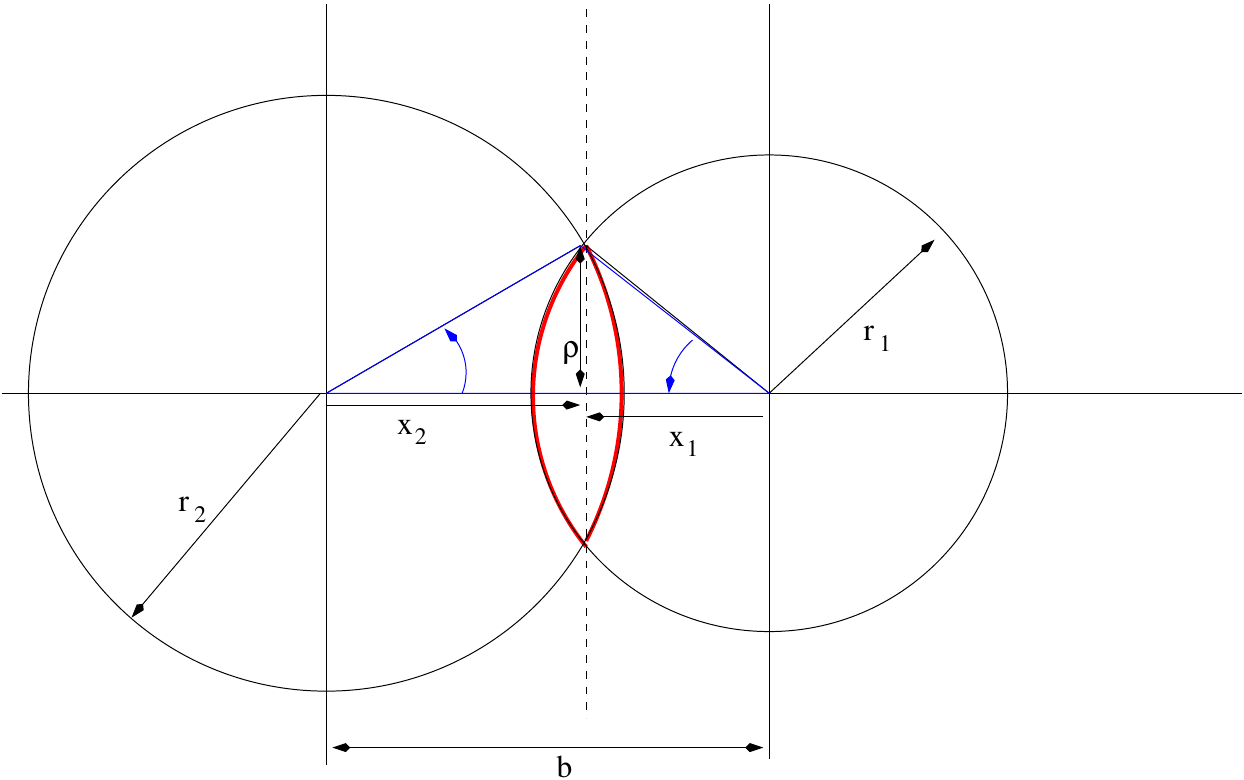}
\caption{
Two intersecting circles of radius $r_1$ and $r_2$ at a distance $b$ with base
radius $\rho$ of the red asymmetric lens.
$b^2> r_1^2+r_2^2$. $x_2>0$. $x_1<0$.
}
\label{fig.circ}
\end{figure}

\begin{figure}
\includegraphics[scale=0.45]{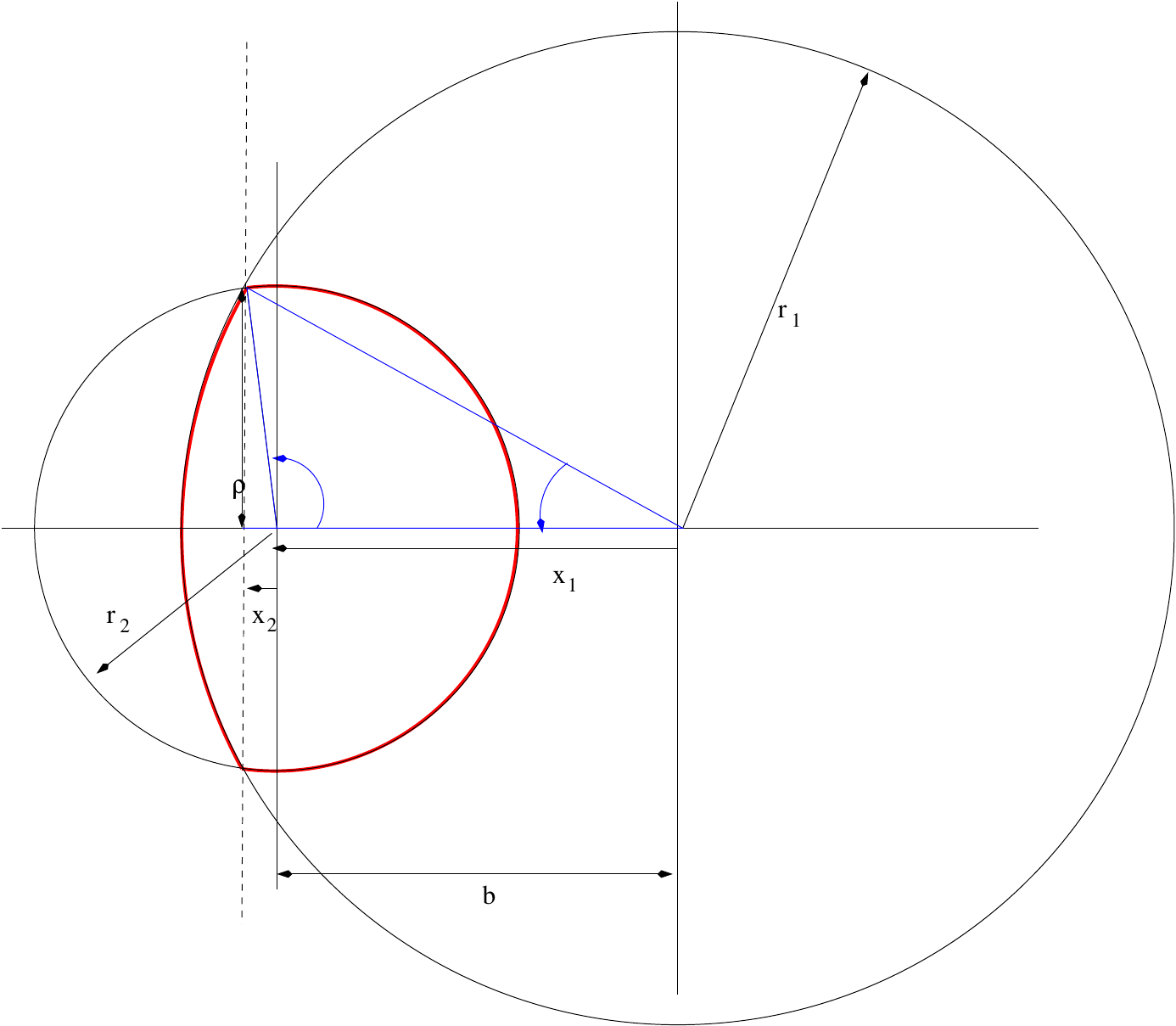}
\caption{
Two intersecting circles of radius $r_1$ and $r_2$ at a distance $b$ with base
radius $\rho$ of the red asymmetric lens.
$b^2< r_1^2+r_2^2$. $x_2<0$. $x_1<0.$
}
\label{fig.circ2}
\end{figure}

The geometry of two planar intersecting circles of radii $r_1$ and $r_2$ at distance $b$
is illustrated in Figures \ref{fig.circ} and \ref{fig.circ2}.
Only the cases with non-vanishing overlap, $b<r_1+r_2$, are of interest here.
We also assume that the circle rims intersect, which means $b+r_2>r_1$ \cite{PetitjeanDG}.

The (right) circle of radius $r_1$ is placed at the center of coordinates:
\begin{equation}
x^2+y^2=r_1^2.
\label{eq.xyr1}
\end{equation}
The (left) circle of radius $r_2$ is placed at $(-b,0)$:
\begin{equation}
(x+b)^2+y^2=r_2^2.
\end{equation}
Solving the first equation for $y^2$, insertion in the second equation and solving for $x$
gives for the horizontal coordinates
of the lens position
\begin{equation}
x_1=\frac{r_2^2-r_1^2-b^2}{2b}<0;\quad x_2=b+x_1 = \frac{r_2^2-r_1^2+b^2}{2b}.
\label{eq.lensOrd}
\end{equation}
$x_2$ is positive in Figure \ref{fig.circ}, negative in Figure \ref{fig.circ2}.
The radius $\rho\ge 0$ of the lens in the right triangle of sides $r_1$, $x_1$ and $\rho$ 
is the associated $y$-value from \eqref{eq.xyr1},
\begin{equation}
\rho=\sqrt{r_1^2-x_1^2}=\sqrt{r_2^2-x_2^2}
=\frac{\sqrt{[(b+r_2)^2-r_1^2][r_1^2-(r_2-b)^2]}}{2b}
.
\label{eq.rho}
\end{equation}
Imaginary
values of $\rho$ are numerical indicators for circles with non-intersecting rims,
because they either are too far apart or one circle lies entirely
within the other.
\begin{rem}
If the value of $x_1$ becomes less than $-b-r_2$ (indicating that the intersection
would be to the left of the left circle), the factor $(b+r_2)^2-r_1^2$ in the
discriminant of this radix becomes negative. If the value of $x_2$ becomes larger
than $r_2$ (indicating that the intersection is to the right of the left circle),
the factor $r_1^2-(r_2-b)^2$ in the discriminant becomes negative. 
Because the polynomial
of the discriminant is a symmetric function of $r_1$ and $r_2$, the
equivalent criteria apply also for the right circle. 
\end{rem}
\begin{rem}
The discriminant of the square root is a polynomial of order 4 in $b$.
Imagine the smaller circle wandering from left to right with decreasing
$b$ in front of the larger circle. There are 4 positions associated
with the roots of this polynomial where 
the two circles have only one point in common and where $\rho$ becomes zero:
2 positions where the circles barely touch, and 2 positions where
the smaller circle is inside the larger circle.
\end{rem}

The angles $\alpha_{1,2}$ at which $\rho$ appears from the centers of the circles
are blue in Figures \ref{fig.circ} \and \ref{fig.circ2}:
\begin{equation}
\sin\alpha_1=\frac{\rho}{r_1};\quad \sin\alpha_2=\frac{\rho}{r_2}.
\end{equation}
\begin{equation}
\cos\alpha_1=\frac{-x_1}{r_1};\quad \cos\alpha_2=\frac{x_2}{r_2}; \quad 0\le \alpha_1,\alpha_2 \le \pi.
\end{equation}
In Figure \ref{fig.circ} $\alpha_2<\pi/2$; 
in Figure \ref{fig.circ2} $\alpha_2>\pi/2$.

The area of the intersection delineated in red
is the sum of the areas of the two circular segments
with radii $r_1$ and $r_2$ \cite[3.76]{Bronstein3}.
Using the principle of inclusion-exclusion, the area of the blue
triangle in Figure \ref{fig.circ}
is the area of the left and right circular sectors minus half the 
area $A$ of the red lens \cite{NelsonApJ174}:
\begin{equation}
\frac12 \rho b = \frac12 \alpha_2 r_2^2+\frac12 \alpha_1 r_1^2 -\frac12 A,
\end{equation}
\begin{equation}
A\left(r_1,r_2,b\right) = \alpha_2 r_2^2+\alpha_1 r_1^2 -\rho b ,
\label{eq.A}
\end{equation}
where $\alpha_{1,2}$ are measured in radians.

A closer inspection of \ref{fig.circ2} shows that the same formula holds.

\section{Elliptic Integral of the Triangular region} \label{app.Atri}
\subsection{4 real roots}
Byrd's reduction of the integral for $v_\Delta$ to standard form is \cite[256.38]{Byrd}
\begin{equation}
\int_b^y\sqrt{(a-t)(t-b)(t-c)(t-d)} dt
=(b-c)^2(a-b)(b-d)\alpha^2 g\int_0^{u_1}\frac{\sn^2u \cn^2u \dn^2 u}{(1-\alpha^2 \sn^2 u)^4}du
\label{eq.iB1}
\end{equation}
where 
\begin{equation}
\alpha^2=(a-b)/(a-c)<1,
\end{equation}
\begin{equation}
k^2=\frac{(a-b)(c-d)}{(a-c)(b-d)},
\end{equation}
\begin{equation}
g=2/\sqrt{(a-c)(b-d)}, 
\end{equation}
\begin{equation}
\varphi = \am u_1=\arcsin\sqrt{\frac{(a-c)(y-b)}{(a-b)(y-c)}},
\label{eq.phiB1}
\end{equation}
\begin{equation}
\sn u_1 = \sin\varphi.
\end{equation}

\cite[362.25]{Byrd}
\begin{multline}
\int_0^{u_1}\frac{\sn^2u \cn^2u \dn^2 u}{(1-\alpha^2 \sn^2 u)^4}du
\\
=\frac{1}{\alpha^6}[-k^2\Pi(u,\alpha^2)+(3k^2-\alpha^2k^2-\alpha^2)V_2
+(2\alpha^2 k^2+2\alpha^2-3k^2-\alpha^4)V_3+(\alpha^2-1)(\alpha^2-k^2)V_4]
\end{multline}
with \cite[336]{Byrd}
\begin{equation}
V_0=F(\varphi,k);
\label{eq.V0}
\end{equation}
\begin{equation}
V_1=\Pi(\varphi,\alpha^2,k);
\end{equation}
\begin{multline}
V_2(u)
=\frac{1}{2(\alpha^2-1)(k^2-\alpha^2)}
\big[
\alpha^2E(u)+(k^2-\alpha^2)u
\\
+(2\alpha^2k^2+2\alpha^2-\alpha^4-3k^2)\Pi(\varphi,\alpha^2,k)-\frac{\alpha^4 \sn u \cn u \dn u}{1-\alpha^2 \sn^2 u}
\big];
\label{eq.Vm2}
\end{multline}
\begin{multline}
V_{m+3}
=\frac{1}{2(m+2)(1-\alpha^2)(k^2-\alpha^2)}
\big[
(2m+1)k^2 V_m
\\
+2(m+1)(\alpha^2k^2+\alpha^2-3k^2)V_{m+1}
\\
+(2m+3)(\alpha^4-2\alpha^2k^2-2\alpha^2+3k^2)V_{m+2}
+\frac{\alpha^4 \sn u \cn u \dn u}{(1-\alpha^2 \sn^2 u)^{m+2}}
\big];
\label{eq.Vm3}
\end{multline}
Note that in the application of this manuscript the upper limit $y$ in \eqref{eq.iB1} equals $a$,
so $\varphi=\pi/2$ in \eqref{eq.phiB1}, $\sn u_1=1$, $\cn u_1=0$, the Elliptic Integrals
are \emph{Complete} Elliptic Integrals, and the Jacobian Elliptic Functions
in \eqref{eq.Vm2}--\eqref{eq.Vm3} do not need to be evaluated.

\subsection{2 complex-conjugated roots}
The case of complex conjugated $d=c^*$ is expanded as
\begin{multline}
\int_b^y \sqrt{(a-t)(t-b)(t-c)(t-c^*)}dt
=
\int_b^y \frac{(a-t)(t-b)(t-c)(t-c^*)}{\sqrt{(a-t)(t-b)(t-c)(t-c^*)}}dt
\\
=
-\int_b^y \frac{t^4}{\sqrt{(a-t)(t-b)(t-c)(t-c^*)}}dt
+(a+b+2\Re c)\int_b^y \frac{t^3}{\sqrt{(a-t)(t-b)(t-c)(t-c^*)}}dt
\\
-(
ab +2a\Re c +2b\Re c +|c|^2
)\int_b^y \frac{t^2}{\sqrt{(a-t)(t-b)(t-c)(t-c^*)}}dt
\\
-(
2ab\Re c +a|c|^2 +b|c|^2
)\int_b^y \frac{t}{\sqrt{(a-t)(t-b)(t-c)(t-c^*)}}dt
-ab|c|^2
\int_b^y \frac{1}{\sqrt{(a-t)(t-b)(t-c)(t-c^*)}}dt
\end{multline}
and relegated to the formulas in App.\ \ref{app.Csec}.

\section{Elliptic Integral of Cone Sector}  \label{app.Csec} 
\subsection{4 real roots}
The integral \eqref{eq.ellcone} is from $t=t_1$ to $t=t_3$, where the argument of the $\arccos$ is $-1$ and $1$.
The following integral needed to be evaluated for $m=2$ and $m=4$ \cite[257.11]{Byrd}:
\begin{equation}
\int_y^a \frac{t^m}{\sqrt{(a-t)(t-b)(t-c)(t-d)}}dt
=a^mg\int_0^{u_1} \frac{(1-\alpha_1^2 \sn^2 u)^m}{(1-\alpha^2 \sn^2 u)^m}du
\label{eq.iB2}
\end{equation}
where
\begin{equation}
\alpha_1^2=\frac{(b-a)d}{a(b-d)}
\end{equation}
\begin{equation}
\alpha^2=\frac{b-a}{b-d}<0
\label{eq.alpha25700}
\end{equation}
\begin{equation}
\varphi = \am u_1=\arcsin\sqrt{\frac{(b-d)(a-y)}{(a-b)(y-d)}},
\label{eq.phiB2}
\end{equation}
The definitions for $k^2$, $g$ and $\sn u_1$ are the same as in Appendix \ref{app.Atri}.
The right-hand side of \eqref{eq.iB2} is \cite[340.04]{Byrd}
\begin{equation}
\int \frac{(1-\alpha_1^2 \sn^2 u)^m}{(1-\alpha^2 \sn^2 u)^m}du
=\frac{\alpha_1^{2m}}{\alpha^{2m}}
\sum_{j=0}^m \binom{m}{j} \frac{(\alpha^2-\alpha_1^2)^j }{\alpha_1^{2j}}V_j,
\label{eq.34004}
\end{equation}
and the $V_j$ given by \eqref{eq.V0}--\eqref{eq.Vm3}.
Note that in the application of this manuscript the lower limit $y$ in \eqref{eq.iB2} equals $b$,
so $\varphi=\pi/2$ in \eqref{eq.phiB2}, $\sn u_1=1$, $\cn u_1=0$, the Elliptic Integrals
are \emph{Complete} Elliptic Integrals, and the Jacobian Elliptic Functions
in \eqref{eq.Vm2}--\eqref{eq.Vm3} do not need to be evaluated.

\subsection{2 complex-conjugated roots}
If $c$ and $d=c^*$ are a pair of complex-conjugate values, 
the applicable entry is \cite[259.03]{Byrd}
\begin{equation}
\int_b^y \frac{t^m}{(a-t)(t-b)(t-c)(t-c^*)}dt
=g\frac{(aB+bA)^m}{(A-B)^m}\sum_{j=0}^m \binom{m}{j}\alpha_2^{m-j}(\alpha-\alpha_2)^j 
\int_0^{u_1} \frac{du}{(1+\alpha \cn u)^j},
\label{eq.25903}
\end{equation}
where $b_1=\Re c$, $a_1=\Im c$, $A^2=(a-b_1)^2+a_1^2$, $B^2=(b-b_1)^2+a_1^2$,
$g=1/\sqrt{AB}$, $\alpha=(A-B)/(A+B)$, $\alpha_2=(bA-aB)/(aB+bA)$, and \cite[341.05]{Byrd}
\begin{equation}
R_m\equiv \int \frac{du}{(1+\alpha \cn u)^m}.
\label{eq.Rm}
\end{equation}
The $R_m$ are recursively
\begin{equation}
R_{-2}=\frac{1}{k^2}[(k^2-\alpha^2k'^2)u+\alpha^2E(u)+2\alpha k\arccos(\dn u)];
\end{equation}
\begin{equation}
R_{-1}=u+\frac{\alpha}{k}\arccos(\dn u);
\end{equation}
\begin{equation}
R_0=u;
\end{equation}
\begin{equation}
R_1=\frac{1}{1-\alpha^2}\left[\Pi(\varphi,\frac{\alpha^2}{\alpha^2-1},k)-\alpha f_1\right]
\end{equation}
at $u=F(\varphi,k)$, $E(u)=E(\varphi,k)$, and \cite[361.54]{Byrd}
\begin{equation}
f_1=\left\{
\begin{array}{ll}
\sqrt{\frac{1-\alpha^2}{k^2+k'^2\alpha^2}}
\arctan\left[\sqrt{\frac{k^2+k'^2\alpha^2}{1-\alpha^2}}\sd u\right], & \alpha^2/(\alpha^2-1)<k^2 ;\\
\sd u, & \alpha^2/(\alpha^2-1)=k^2 ;\\
\frac12\sqrt{\frac{\alpha^2-1}{k^2+k'^2\alpha^2}}\ln
\frac{\sqrt{k^2+k'^2\alpha^2}\dn u+\sqrt{\alpha^2-1}\sn u}{\sqrt{k^2+k'^2\alpha^2}\dn u-\sqrt{\alpha^2-1}\sn u},
& \alpha^2/(\alpha^2-1)>k^2 ;\\
\end{array}
\right.
\end{equation}
\cite[341.05]{Byrd}
\begin{multline}
R_m=\frac{1}{(m-1)(\alpha^2-1)(k^2+\alpha^2k'^2}\Big\{
(3-2m)[\alpha^2(1-2k^2)+2k^2]R_{m-1}
\\
+2(5-2m)(k^2R_{m-3}
+(m-1)(6k^2+\alpha^2-2k^2\alpha^2)R_{m-2}
+(m-3)k^2R_{m-4}+\frac{\alpha^3\sn u\dn u}{(1+\alpha \cn u)^{m-1}}
\Big\}
.
\end{multline}

\section{Elliptic Integral of Sphere Sector}  \label{app.Ssec} 
\subsection{4 real roots}
Eq. \ref{eq.d} requires the integrals
\cite[3.151.7]{GR}
\cite[257.39]{Byrd}
\begin{equation}
\int_y^a \frac{dt}{(p-t)^m\sqrt{(a-t)(t-b)(t-c)(t-d)}}
=\frac{g}{(p-a)^m}\int_0^{u_1}\frac{(1-\alpha^2 \sn^2u)^m}{(1-\alpha_1^2 \sn^2u)^m}du
\end{equation}
where $p\neq a$, $\alpha_1^2\equiv (p-d)(a-b)/(a-p)/(b-d)$,
$\alpha^2$ defined in \eqref{eq.alpha25700}, and the 
integral of the right hand side reduced in \eqref{eq.34004}.

\subsection{2 complex-conjugated roots}
If $d=c^*$ are complex conjugated in the previous integral \cite[259.04]{Byrd}
\begin{multline}
\int_b^y \frac{1}{(t-p)^m\sqrt{(a-t)(t-b)(t-c)(t-c^*)}}dt
\\
=\frac{(A+B)^m}{[A(b-p)-B(a-p)]^m}g\sum_{j=0}^m \binom{m}{j}
\alpha_1^{m-j}(\alpha-\alpha_1)^j
\int_0^{u_1} 
\frac{du}{(1+\alpha \cn u)^j}
\end{multline}
where $A$ and $B$ are defined after
\eqref{eq.25903}, $\alpha\equiv (bA-aB+pB-pA)/(aB+bA-pA-pB)$, $\alpha_1\equiv (A-B)/(A+B)$,
and the right hand side is evaluated with \eqref{eq.Rm}.

\subsection{C++ implementation}
The computation of the volume is implemented in the C++ source code in
the ancillary directory, using the GNU scientific library  (GSL) to evaluate
the elliptic integrals \cite{GSL}. For Linuxes the minimum package names
depend on the distribution. For openSUSE compiler and GSL are retrieved
with \texttt{zypper install gcc-c++ cpp gsl-devel make}, for example, 
on Ubuntu like \texttt{apt install g++ cpp libgsl-dev make}.
The \texttt{Makefile} compiles two main programs, \texttt{sphereCylVol}
and \texttt{sphereConeVol}: 
\begin{itemize}
\item
The volume of the intersection of a
sphere and a cylinder 
\cite{LamarcheCPC59}
is calculated with the call

\small
\begin{verbatim}
sphereCylVol r R b
\end{verbatim}
\normalsize

with three floating point parameters: $r$ is the sphere radius, $R$ the cylinder radius, and $b$ the impact parameter.
Viviani's volume for example
is computed with \texttt{sphereCylVol 1. 0.5 0.5}.
\item
The volume of intersection of a sphere and a cone is calculated with one of
 \small

\begin{verb}
sphereConeVol [-N samples [-v]] [-r radius] [-p phiDegrees] sx sy sz ax ay az dx dy dz
\end{verb}

\begin{verb}
sphereConeVol [-N samples [-v]] [-r radius] [-p phiDegrees] ax ay az dx dy dz
\end{verb}

\begin{verb}
sphereConeVol [-N samples [-v]] [-r radius] [-p phiDegrees] ax ay az
\end{verb}

\normalsize

The option \texttt{-r} followed by a positive floating-point number
specifies the sphere radius $R$. If the argument is not used the radius
is assumed to be 1.

The option \texttt{-p} followed by a positive floating-point number
specicies the cone half angle $\varphi$ in degrees. If the argument is not
used it is assumed to be 45.

The 9, 6 or 3 trailing arguments are signed floating point numbers with groups of Cartesian
x, y and z coordinates. (If at least one of the numbers is negative, a double-dash \texttt{--}
should be inserted after the options to disambiguate the meaning of their minus-sign
and the dashes of the options.) The triple \texttt{sx sy sz} are the Cartesian coordinates
of the sphere center. If the triple is missing, $(0,0,0)$ is assumed. The triple \texttt{ax ay az}
are the Cartesian coordinates of the cone apex. The triple \texttt{dx dy dz} are the Cartesian
coordinates of the direction of the cone axis. (The length of that vector does not
need to be normalized to unity but must be nonzero.) If \texttt{dx dy dz} are absent, 
the direction $(0,0,1)$ parallel to the z-axis is assumed.

The option \texttt{-N} is a debugging option and lets the program compute the 
approximate (!) volume by slicing the sphere into that many pieces and adding the areas of the 
circular intersections with the cone in the sense of a Simpson summation.
The higher the
integer argument \textit{samples}, the more accurate the result. If the argument is
not used, the integrals of this manuscript are evaluated. For increasingly large \textit{samples}
both results ought to converge.

If the samples argument is a negative number, the approximate area (!)
of the sphere surface
is calculated which is inside the cone. 
[There is no equivalent analytical
evaluation of areal integrals in this manuscript.]
This samples
the front surface and also the back surface. 
In satellite imaging the back surface would not be visible;
if the option \texttt{-v} is also used, only the visible area
is accumulated (i.e. the parts where the vector from sphere center to
surface and the vector from apex to sphere surface have a dot product less than zero).

\item
The the validity of the analytical integrals is investigated with

\begin{verb}
sphereConeVol -t [-N samples]
\end{verb}

which runs a triple loop over a finite grid over the three parameters $\hat d$, $\hat b$ and $\varphi$
and prints for each point $\hat b$, $\hat d$, $\varphi$ (in radians),
the result obtained
by the elliptic integrals and the result obtained by a Simpson rule.
If the two values of the volume differ by more than $10^{-6}$,
an additional exclamation mark and the difference is printed.

The numerical Simpson rule slices
the sphere in \textit{samples} pieces. If the \texttt{-N} option is not used, a
value of 100,000 is assumed. 

This is equivalent to a few hundred separate calls of \texttt{sphereConeVol} with and
without the \texttt{-N} option.
\end{itemize}

\bibliographystyle{amsplain}
\bibliography{all}

\providecommand{\bysame}{\leavevmode\hbox to3em{\hrulefill}\thinspace}
\providecommand{\MR}{\relax\ifhmode\unskip\space\fi MR }
\providecommand{\MRhref}[2]{%
  \href{http://www.ams.org/mathscinet-getitem?mr=#1}{#2}
}
\providecommand{\href}[2]{#2}
\begin{thebibliography}{10}

\bibitem{GSL}
\emph{{GNU} scientific library}, 2023, http://www.gnu.org/software/gsl/.

\bibitem{AS}
Milton Abramowitz and Irene~A. Stegun (eds.), \emph{Handbook of mathematical
  functions}, 9th ed., Dover Publications, New York, 1972. \MR{0167642}

\bibitem{BoersmaPKAW64}
J.~Boersma and W.~Kamminga, \emph{Calculation of the volume of intersection of
  a sphere and a cylinder}, Proc. Konink. Akad. Wetensch. A \textbf{64} (1961),
  496--507. \MR{0133036}

\bibitem{Bronstein2}
I.~N. Bronstein and K.~A. Semendjajew, \emph{Teubner's {T}aschenbuch der
  {M}athematik}, Teubner, 1996. \MR{1446529}

\bibitem{Bronstein3}
I.~N. Bronstein, K.~A. Semendjajew, G.~Musio, and H.~Muehlig, \emph{Handbook of
  mathematics}, Springer-Verlag, 2007.

\bibitem{Byrd}
Paul~F. Byrd and Morris~D. Friedman, \emph{Handbook of elliptical integrals for
  engineers and physicists}, 2nd ed., Die Grundlehren der mathematischen
  Wissenschaften in Einzeldarstellungen, vol.~67, Springer, Berlin,
  G\"ottingen, 1971. \MR{0277773}

\bibitem{GR}
I.~Gradstein and I.~Ryshik, \emph{Summen-, {P}rodukt- und {I}ntegraltafeln},
  1st ed., Harri Deutsch, Thun, 1981. \MR{0671418}

\bibitem{LamarcheCPC59}
Fran\c{c}ois Lamarche and Claude Leroy, \emph{Evaluation of the volume of
  intersection of a sphere with a cylinder by elliptic integrals}, Comp. Phys.
  Commun. \textbf{59} (1990), no.~2, 359--369. \MR{1058694}

\bibitem{MillerACMTG6}
James~R. Miller, \emph{Geometric approaches to nonplanar quadric surface
  intersection curves}, ACM Trans. Graphics \textbf{6} (1987), no.~4, 274--307.

\bibitem{NelsonApJ174}
Burt Nelson and Walter~D. Davis, \emph{Eclipsing-binary solutions by sequential
  optimization of the parameters}, Astroph. J. \textbf{174} (1972), 617--628.

\bibitem{PetitjeanDG}
Michel Petitjean, \emph{Spheres unions and intersections and some of their
  applications in molecular modeling}, pp.~61--83, Springer, 2013. \MR{3051948}

\bibitem{RuffJAM10}
I.~Ruff, \emph{The intersection of a cone and a sphere: A contribution to the
  geometry of satellite viewing}, J. Appl. Meteor. \textbf{10} (1971), no.~3,
  607--609.

\bibitem{SarragaCVGI22}
Ramon~F. Sarraga, \emph{Algebraic methods for intersections of quadric surfaces
  in {GMSOLID}}, Computer Vis., Graphics Image Proc. \textbf{22} (1983), no.~2,
  222--238.

\bibitem{SchmidCGS20}
Erwin Schmid, \emph{The earth as viewed from a satellite}, National Geodetic
  Survey \textbf{20} (1962), 1--21.

\bibitem{ShahSIAMR30}
G.~A. Shah, \emph{Problem 88-10, volume of intersection of a cone with a
  sphere}, SIAM Rev. \textbf{30} (1988), no.~2, 310--311.

\end{thebibliography}

\end{document}